\documentclass[12pt]{amsart}
\usepackage{amsmath,amsthm,amssymb,amsfonts,amscd}
\usepackage{mathrsfs}
\usepackage{bbm}
\usepackage{bbding}
\usepackage{graphicx,latexsym}
\usepackage[backref=page]{hyperref}
\usepackage{hyperref}
\usepackage{geometry}
\usepackage{color}
\usepackage{xcolor}
\usepackage{picture,epic}
\usepackage{tikz}
\usepackage{verbatim}
\usepackage{}
\usepackage{enumitem}
\usepackage[T1]{fontenc}

\definecolor{rev1color}{rgb}{1, 0, 0}  
\definecolor{rev2color}{rgb}{1, 0, 0}  

\newif\ifroundone
\roundonefalse 

\newcommand{\revised}[1]{%
	\ifroundone
	\textcolor{rev1color}{#1}%
	\else
	#1%
	\fi
}

\newif\ifroundtwo
\roundtwofalse 

\newcommand{\revisedtwo}[1]{%
	\ifroundtwo
	\textcolor{rev2color}{#1}%
	\else
	#1%
	\fi
}

\numberwithin{equation}{section}

\theoremstyle{plain}
\newtheorem{theorem}{Theorem}[section]
\newtheorem{lemma}[theorem]{Lemma}
\newtheorem{corollary}[theorem]{Corollary}
\newtheorem{proposition}[theorem]{Proposition}

\theoremstyle{definition}

\theoremstyle{remark}
\newtheorem{remark}[theorem]{Remark}

\renewcommand{\Re}{\operatorname{Re}}
\renewcommand{\Im}{\operatorname{Im}}

\newcommand{\GL}{\operatorname{GL}}

\renewcommand{\mod}{\operatorname{mod}\ }

\newcommand{\dd}{\mathrm{d}}

\makeatletter
\@namedef{subjclassname@2020}{%
	\textup{2020} Mathematics Subject Classification}
\makeatother

\begin{document}	
	
	
	\baselineskip=17pt
	
	\title{\revised{On the Second Moment of Twisted Higher Degree $L$-functions}}
	
	\author[H. Gou]{Haozhe Gou}
	\address{School of Mathematics\\ Shandong University \\ Jinan \\ Shandong 250100 \\China}
	\email{hodgegou@mail.sdu.edu.cn}
	
	\author[L. Li]{Liangxun Li}
	\address{Data Science Institute and School of Mathematics \\ Shandong University \\ Jinan \\ Shandong 250100 \\China}
	\email{lxli@mail.sdu.edu.cn}	

	\begin{abstract}
   \revised{Assuming the Ramanujan conjecture, the zero density estimate and \revisedtwo{some  subconvexity type bound}}, we  describe a general method to obtain the log-saving upper bound for the second moment of standard twisted \revised{higher degree} $L$-function in the $q$-aspect.
    Specifically, let $L(s, F)$ be a standard $L$-function of degree $d\geq3$\revisedtwo{. \revised{Under} these foundational hypotheses,} the bound
		\[
		\sideset{}{^*}{\sum}_{{\chi (\mod q)}}\Big|L\big(\frac{1}{2}, F\times \chi \big)\Big |^2\ll_{F,\eta} \frac{q^{\frac{d}{2}}}{\log^{\eta}q}
		\]
		holds for some small $\eta>0$.
	\end{abstract}
	
\maketitle

\section{Introduction}\label{sec:Intr}
A fundamental problem in number theory is estimating the values of \(L\)-functions at the central point.  Many \(L\)-functions are organized into families, and their moments often encode some deep properties of the families.
In recent years, the families of twisted $L$-functions have attracted great attention.
There are many conjectures concerning the asymptotic behaviour for the moments of twisted $L$-functions.
And in the case of lower degree, \revisedtwo{some of them have been settled by using advanced techniques.} These results also produce various important applications such as the non-vanishing problem, extreme values of $L$-functions, etc.

Let $q\geq 2$ be an integer and $\pi$ be an irreducible automorphic representation on $\GL_d$ ($d\geq 1$) over $\mathbb{Q}$ with unitary central character. The $L$-function of $L(s, \pi)$ twisted by a primitive Dirichlet character $\chi$ modulo $q$ is
\begin{equation*}
L(s, \pi\times \chi)=\sum_{n\geq 1}\frac{\lambda_\pi(n)\chi(n)}{n^s}, \quad \Re(s)>1,
\end{equation*}
where $\lambda_{\pi}(n)$ is the $n$-th Dirichlet coefficient of $L(s,\pi)$.
There is an important problem which is to understand the even moments of the above $L$-function at the central value when $\chi$ run through all the primitive Dirichlet characters modulo $q$:
\begin{equation}\label{eqn: 2kmoment}
\sideset{}{^*}\sum_{\chi(\mod q)}\Big|L(\frac{1}{2},\pi\times \chi)\Big|^{2k},
\end{equation}
\revisedtwo{where $k \ge 1$ is a positive integer.}
When the degree $d=1$ and $k=1, 2$,  we face with the moments of Dirichlet $L$-functions
\[
\sideset{}{^*}\sum_{\chi(\mod q)}|L(1/2,\chi)|^{2}, \quad\text{and}\quad
\sideset{}{^*}\sum_{\chi(\mod q)}|L(1/2,\chi)|^{4}.
\]
The asymptotic formula for the second moment of Dirichlet $L$-functions has been first proven by Paley \cite{Paley1931}.
In \cite{HB81}, Heath-Brown has computed a similar moment, but for all characters modulo $q$.
If $q$ is a prime number, there are several improvements such as \cite{Conrey07}, \cite{Young11}.
This asymptotic is significant due to the fact that it is close to the non-vanishing of $L(\frac{1}{2}, \chi)$.
By computing the first and the second twisted moment, many positive propositional results for non-vanishing problem have been established, for example, see \cite{BM92}, \cite{IS99}, \cite{KMN22}, etc.
The asymptotic formula of the fourth power moment has been obtained by Heath-Brown \cite{He81}, for $q$ prime, and  by Soundararajan \cite{So05} for general $q.$ Young \cite{Yo11} established an asymptotic formula with power saving for prime moduli $q:$
$$\sideset{}{^*}\sum_{\chi(\mod q)}|L(1/2,\chi)|^4=qP_4(\log q)+O\left(q^{1-\delta +\varepsilon}\right),$$
\revisedtwo{where $P_4$ is a polynomial of degree four with leading coefficient $1/(2\pi^2)$. The admissible exponent} was initially established as 
$\delta=\frac{1}{80}+\frac{\theta}{40}$ with $\theta\leqslant7/64$ 
being an admissible exponent towards the Ramanujan-Petersson conjecture. 
\revised{Subsequent improvements first achieved $\delta=\frac{1}{32}$ in \cite{Bl17}, 
	then further advanced to $\delta=\frac{1}{20}$ by Blomer et al. \cite{BlomerSteklov17}, 
	with the current state-of-the-art $\delta=\frac{1}{14}$ established by Wu \cite{WuMathAnn23}.}

When $d=2$, $k=1$, it leads to consider
\begin{equation}\label{eqn:1.2}
	\sideset{}{^*}\sum_{\chi(\mod q)}|L(1/2,f\times\chi)|^2,
\end{equation}
where $f$ is a fixed $\rm GL_2$ Hecke--Maass cusp form.
Also if we relax $f$ to be a (derivative of an) Eisenstein series, the corresponding formula exactly equals to  the fourth moment of Dirichlet $L$-functions. The second moment \eqref{eqn:1.2} was studied by Stefanicki \cite{St96}, who proved an asymptotic formula with an error term that saves a small power of $\log q$, provided $q$ has only few prime divisors. An asymptotic formula with a $\log \log q$-saving was established by Gao, Khan and Ricotta \cite{G-K-R09} for almost all integers $q$. In \cite{Bl15}, assuming the \revisedtwo{factorable} property of the moduli $q$, Blomer and Mili\'{c}evi\'{c} established the asymptotic with a power-saving error term.
For a prime moduli $q$, Blomer, Fouvry, Kowalski, Michel, Mili\'{c}evi\'{c}, and Sawin successfully obtained a small power-saving error term in their combination works \cite{Bl17,Ko17}.
In particular, by applying the ``$+ab$ shift'' trick and leveraging deep results from algebraic geometry, new bounds for bilinear forms in hyper-Kloosterman sums that go beyond the Fourier-theoretic range are established in \cite{Ko17} (building on previous work). As an application, they \revisedtwo{provide}
$$
\sideset{}{^*}\sum_{\chi(\mod q)}|L(1/2,f\times\chi)|^2=qP_{f}(\log q)+O_{f}(q^{1-\frac{1}{144}+\varepsilon}),
$$
where $P_f(X)$ is a polynomial of degree 1 depending on $f$ only.
These asymptotic results mentioned above also have general versions in the mixed moments (which we will not discuss here).


For the upper bound in the case of higher degrees ($d\geq 3$) and higher moments, by using the approximate functional equation for $L(\frac{1}{2},\pi\times \chi)$ and the orthogonality of Dirichlet characters, we find that
\[
\sideset{}{^*}\sum_{\chi(\mod q)}\Big|L(\frac{1}{2},\pi\times \chi)\Big|^{2k}\ll_{\pi} q^{\frac{dk}{2}+\varepsilon}.
\]
This yields the convexity bound $L(\frac{1}{2},\pi\times \chi)\ll_{\pi} q^{\frac{d}{4}+\varepsilon}.$
Here the analytic conductor $\revised{C(\pi\times\chi)\asymp}_\pi q^d.$
This type bound sometimes is called the trivial bound for $L$-functions. Assuming the generalized Lindel\"of hypothesis in the conductor aspect $L(\frac{1}{2},\pi\times \chi)\ll q^{\varepsilon},$ which is a consequence of the generalized Riemann hypothesis, the upper bound of the moment \eqref{eqn: 2kmoment} is $O(q^{1+\varepsilon}).$  When $d\geq 3$, it is still widely open to seek the non-trivial upper bound for \eqref{eqn: 2kmoment} than $O(q^{\frac{dk}{2}+\varepsilon}).$



In this paper, we are interested in the upper bound for the second moment of more general twisted $L$-function in the $q$-aspect. We describe an axiomatic framework for the $L$-functions considered in this paper.
Let $L(s, F)$
\footnote{
As an example, we can simply recognize $F=\pi$  where $\pi$ is an irreducible automorphic representation on $\GL_d$ over $\mathbb{Q}$ with unitary central character. And $L(s, F)=L(s, \pi)$ is the automorphic $L$-function associated with $\pi$. All of these basic properties are satisfied for $L(s, \pi)$ and its twist $L(s, \pi \times \chi).$ This case will be discussed in Section \S \ref{sec:pi_case}.
}
 be given by the Dirichlet series and Euler product of degree $d\geq 3$,
\begin{equation}\label{eqn: L-def}
L(s, F)=\sum_{n\geq 1}\frac{\lambda_F(n)}{n^s}=\prod_{p}\prod_{j=1}^{d}\Big(1-\frac{\alpha_{j, F}(p)}{p^s}\Big)^{-1},
\end{equation}
and we assume that both of them are absolutely convergent in $\Re(s)>1.$
Here $\lambda_F$ is a multiplicative function defined by
\begin{equation*}
	\lambda_F(p^r)=\sum_{r_1,\cdots, r_d\geq 0\atop r_1+\cdots+r_d=r}\prod_{1\leq j\leq d}\alpha_{j, F}(p)^{r_j}
\end{equation*}
at any prime power $p^r$.
Let $q\geq 2$ and $\chi$ is a primitive Dirichlet character modulo $q$. \revisedtwo{The $L$-function associated to $F$ and twisted by $\chi$ is defined by}
\begin{equation}\label{eqn:L F twist chi}
		L(s, F\times \chi)=\sum_{n\geq 1}\frac{\lambda_F(n)\chi(n)}{n^s}=\prod_{p}\prod_{j=1}^{d}\Big(1-\frac{\alpha_{j, F}(p)\chi(p)}{p^s}\Big)^{-1},
\end{equation}
and we also suppose that they are absolutely convergent in $\Re(s)>1.$ We set
\begin{equation}\label{eqn: L_infty}
	L_{\infty}(s, F\times \chi)=N^{\frac{s}{2}}\prod_{j=1}^{d}\Gamma_{\mathbb{R}}(s+\mu_j),
\end{equation}
where $\Gamma_{\mathbb{R}}(s)=\pi^{-\frac{s}{2}}\Gamma(\frac{s}{2})$, $N$ denotes the conductor which depends only on $F$ and $q$,
the parameters $\{\mu_j\}_{1\leq j\leq d}$ are complex numbers depending on $F$ and $\chi(-1)$ such that
\begin{equation}\label{eqn: theta_f}
	\Re(\mu_j)\geq -1+\theta_F,
\end{equation}
for some fixed positive $\theta_F$.  We assume that $\theta_F$ is in the following range
\begin{equation}\label{eqn: theta_range}
	\frac{1}{2d}<\theta_F<\frac{1}{2}.
\end{equation}
The complete $L$-function is defined by
\begin{equation*}
	\Lambda(s, F\times \chi):=L(s, F\times \chi)L_{\infty}(s, F\times \chi).
\end{equation*}
 \revised{The following properties are assumed for $\Lambda(s, F \times \chi)$:}
	\begin{enumerate}
		\item \revised{(Entirety)} $\Lambda(s, F \times \chi)$ extends to an entire function of order 1;
		\item \revised{(Functional Equation)} For some root number $\kappa = \kappa_{F, \chi}$ with $|\kappa| = 1$,
		\begin{equation}\label{eqn: FE}
			\Lambda(s, F\times \chi)=\kappa \Lambda(1-s, \tilde{F}\times \overline{\chi}),
		\end{equation}
		where $\tilde{F}$ is the dual form satisfying $\lambda_{\tilde{F}}(n) = \overline{\lambda_F(n)}$;
		\item \revised{(Analytic Conductor)} The analytic conductor is given by
		\[
		\revised{C(s,F\times\chi) = N\prod_{j=1}^{d}\left(|s+\mu_j+3|\right).}
		\]
	\end{enumerate}

\revised{These axioms hold for automorphic $L$-functions of $\GL(n)$; see \cite[Section 5.2]{I-K04}.}
In our paper, we consider the central value of $L$-function, that is $s=\frac{1}{2}$. We focus on the conductor aspect for $L(s, F\times \chi)$, so that $\revised{C(F\times\chi):=C(\frac{1}{2},F\times\chi)}\ll_{F} N$. Morever we suppose that the bound
\begin{equation}\label{eqn: N-bound}
	\revised{C(F\times\chi)}\ll_{F} N\ll_{F}  q^{d}
\end{equation}
holds for $q\geq 2$ large.\par
\subsection{Main Results}\label{subsec:mresult}
\revisedtwo{The above properties in the beginning of introduction are satisfied for most standard $L$-functions} in analytic number theory.
For the sake of completeness, we assume that $F$ and its associated $L$-functions satisfy the above properties throughout the paper.
In order to describe our result, we need some extra profound assumptions on $F$ and $L(s, F\times \chi)$.
For the multiplicative function $\lambda_F$, we expect the Ramanujan bound:
\[
|\lambda_F(n)|\leq \tau_d(n),
\leqno(\emph{A1})
\]
for all integer $n\geq 1$, where $\tau_d$ is the $d$-fold divisor function.
This can be induced from the generalized Ramanujan conjecture (GRC) on $L(s, F)$ which asserts that
$|\alpha_{j, F}(p)|\leq 1$ for all $1\leq j\leq d.$
Moreover, we assume that the bound on average:
\[
\sum_{X\leq n\leq X+Y}|\lambda_F(n)|^2\ll_{F} Y,
\leqno(\emph{A2})
\]
\revisedtwo{holds for any $X\geq 2$ large} and $X^{1-\varepsilon}\leq Y\leq X$ for \revised{some fixed
	(but arbitrary) $\varepsilon>0.$}
Also, we require a sieve bound on the mean value:
\[
\sum_{n\leq X\atop (n, \prod_{P_1\leq p\leq P_2}p)=1}|\lambda_F(n)|^2\ll_{F} X\frac{\log P_1}{\log P_2},
\leqno(\emph{A3})
\]
for any $2\leq P_1\leq P_2\leq X.$
The last assumption concerns the zero free region for the $L$-function $L(s, F\times\chi)$.
Let $\mathcal{F}_q$ be a subset (\revisedtwo{including the complete set, may depend on $F$}) of all primitive Dirichlet character modulo $q$,
we expect that  each of $L$-function $L(s, F\times\chi)$ with $\chi\in \mathcal{F}_q$ has the following zero free region:
\[
\prod_{\chi\in \mathcal{F}_q}L(\sigma+it, F\times \chi)\neq 0, \quad \text{ for } \sigma>1-\frac{c_F}{\log^{\beta}(qT)},\ \ {|t|\leq T},
\leqno(\emph{A4})
\]
for any $ \log^{100} q\ll T\ll q^{o(1)}$, where $c_F>0$ is a fixed small constant depending on $F$ and  $0\leq \beta<1$.

\begin{theorem}\label{thm:2ML}
	\revisedtwo{Let $d\geq 3$ and $q\geq 2$. Assume  (\emph{A1}),  (\emph{A2}),  (\emph{A3}) and  (\emph{A4}). Then we have}
	\begin{equation}
		\sideset{}{}{\sum}_{{\chi \in \mathcal{F}_q}}\Big|L\big(\frac{1}{2}, F\times \chi \big)\Big |^2\ll_{F} \frac{q^{\frac{d}{2}}}{\log^{\eta}q}
	\end{equation}
	for any $\eta<1-\beta$.
\end{theorem}
\begin{theorem}\label{thm:2ML2}
	\revisedtwo{Let $d\geq 5$ and $Q\geq 2$. Assume (\emph{A1}),  (\emph{A2}),  (\emph{A3}) and  (\emph{A4}) for all $Q\leq q\leq 2Q$. Then we have}
	\begin{equation}
		\sum_{Q\leq q\leq 2Q}\sideset{}{}{\sum}_{\chi \in \mathcal{F}_q }\Big|L\big(\frac{1}{2}, F\times \chi \big)\Big |^2\ll_{F} \frac{Q^{\frac{d}{2}}}{\log^{\eta}Q}
	\end{equation}
	for any $\eta<1-\beta$.
\end{theorem}

As a straightforward consequence of Theorem \ref{thm:2ML}, where \(\mathcal{F}_q\) is chosen to be the complete set of all primitive Dirichlet characters modulo \(q\), it yields the  weak subconvexity bound. For the general case of weak subconvexity and subconvexity problem, we refer to ~\cite{Sound10}, \cite{S-T19} and \cite{Nelson21}. 

\begin{corollary}\label{cor:weaksb}
	Let the assumption be the same as in Theorem \ref{thm:2ML}, and (A4) holds for all primitive Dirichlet  characters modulo $q$, then we have
	\begin{equation}
		L\big(\frac{1}{2}, F\times \chi \big)\ll_{F} \frac{q^{\frac{d}{4}}}{\log^{\eta}q}
	\end{equation}
	for any $\eta<\frac{1-\beta}{2}$.
\end{corollary}

\revisedtwo{As we all know,} establishing  (\emph{A4}) with $\beta<1$ is a standing challenge in the theory of $L$-function. By the classic analytic approach, $\beta\geq 1$ may hold in the zero free region (\emph{A4}). And $\beta$ can be $\frac{2}{3}+o(1)$ in the case of Riemann $\zeta$-function on $t$-aspect. In order to achieve a feasible method to establish the above theorems (under GRC), we assume the following two conditions to avoid (\emph{A4}). For $1-\frac{1}{2024}<\sigma<1$, $T\geq 1$, let $\mathcal{N}_F(\sigma, T; \chi)$ be the number of zeros of $L(s, F\times \chi)$ in the region $\Re(s)\geq \sigma$ and $|\Im(s)|\leq T.$
The first condition is a zero density estimate: for all $1\leq T\leq q$,
\[
\sideset{}{^*}{\sum}_{{\chi (\mod q)}}\mathcal{N}_F(\sigma, T; \chi)\ll_{F} q^{C_F(1-\sigma)}\log^A q,
\leqno(\emph{A5})
\]
where $C_F$ is a large fixed positive constant depending on $F$ and $A$ is an arbitrarily large constant.
If we consider the average over $Q\leq q\leq 2Q$,  then we need a stronger version of (\emph{A5}):
for all $1\leq T\leq Q$, we have
\[
\sum_{Q\leq q\leq 2Q}\sideset{}{^*}{\sum}_{{\chi (\mod q)}}\mathcal{N}_F(\sigma, T; \chi)\ll_{F} Q^{C_F'(1-\sigma)}\log^A Q
\leqno(\emph{A5'})
\]
for some large fixed constant $C_F'>0$. \revised{Indeed, (\emph{A5}) and (\emph{A5'}) is known for us via the work of Lemke Oliver and Thorner \cite{L-T19}.}

The second condition is a \revised{``weak''} type subconvexity bound for $L(\frac{1}{2}, F\times \chi)$ on $q$-aspect:
\[
L\Big(\frac{1}{2}, F\times \chi\Big)\ll_{F} \frac{q^{\frac{d}{4}}}{\exp(\log^\delta q)}
\leqno(\emph{A6})
\]
for all Dirichlet primitive character $\chi$ modulo $q$,
where $\delta\in (0, 1)$ is a constant. \revised{While labeled as ``weak'', the subconvexity bound (A6) is exponentially stronger than the Soundararajan--Thorner bound \cite{S-T19}, though it is indeed weaker than  the power-saving subconvexity bound.} 


\begin{theorem}\label{thm:2ML'}
Let $d\geq 3$ and $q\geq 2$, assuming (A1), (A2), (A3), (A5) and (A6). Then we have
	\begin{equation}
		\sideset{}{^*}{\sum}_{{\chi (\mod q)}}\Big|L\big(\frac{1}{2}, F\times \chi \big)\Big |^2\ll_{F} \frac{q^{\frac{d}{2}}}{\log^{\eta}q}
	\end{equation}
	for any $\eta< \delta$.
\end{theorem}
\begin{theorem}\label{thm:2ML2'}
Let $d\geq 5$ and $Q\geq 2$, assuming (A1), (A2), (A3), (A5') and (A6) holds for all $Q\leq q\leq 2Q$. Then we have
	\begin{equation}
	\sum_{Q\leq q\leq 2Q}\sideset{}{^*}{\sum}_{\chi (\mod q)}\Big|L\big(\frac{1}{2}, F\times \chi \big)\Big |^2\ll_{F} \frac{Q^{\frac{d}{2}}}{\log^{\eta}Q}
	\end{equation}
	for any $\eta< \delta$.
\end{theorem}
We will demonstrate in Section \S\ref{sec: a5a6 to a4} how to replace  (\emph{A4}) with  (\emph{A5}) and  (\emph{A6}), thereby proving Theorems \ref{thm:2ML'} and \ref{thm:2ML2'}. As an application, we have the following theorem.

\begin{theorem}\label{thm:pi_case}
 Let $\pi$ be an irreducible automorphic representation on $\GL_d$ with $d\geq 3$.  Assuming GRC on $\pi$ and the subconvexity bound
 \[
 L\Big(\frac{1}{2}, \pi\times \chi\Big)\ll_{\pi} q^{\frac{d}{4}}e^{-\log^{\delta}q}
 \]
 for all primitive character $\chi$ modulo $q$ where $\delta \in (0,1)$.
 Then we have
  \begin{equation}
		\sideset{}{^*}{\sum}_{{\chi (\mod q)}}\Big|L\big(\frac{1}{2}, \pi\times \chi \big)\Big |^2\ll_{\pi} \frac{q^{\frac{d}{2}}}{\log^{\eta}q}
	\end{equation}
for any $\eta<\delta$. Moreover, with the same assumptions for all $Q\leq q\leq 2Q$ but $d\geq 5$, we have
 \begin{equation}
		\sum_{Q\leq q\leq 2Q}\sideset{}{^*}{\sum}_{{\chi (\mod q)}}\Big|L\big(\frac{1}{2}, \pi\times \chi \big)\Big |^2\ll_{\pi} \frac{Q^{\frac{d}{2}}}{\log^{\eta}Q}
	\end{equation}
for any $\eta<\delta$.
\end{theorem}
This is a special case of $F = \pi$. We will verify it in Section~\S\ref{sec:pi_case}.\par
\subsection{Second Moment Estimates}\label{sec: 2me}
The key ingredient in the proof of the above theorems is based on the second moment estimate on twisted Dirichlet coefficients.
\begin{proposition}\label{prop:2m}
  Let $q\geq 2$. Assuming (A1), (A2), (A3) and (A4),  then for $X\geq q^{1+{\varepsilon}}$ with $\log X\asymp \log q$, we have
  \begin{equation}\label{eqn: lambdaF2mq}
	\sideset{}{}{\sum}_{{\chi \in\mathcal{F}_q}}\Big|\sum_{n\leq X}\lambda_F(n)\chi(n)\Big|^2\ll_F \frac{X^2}{\log^\eta X},
  \end{equation}
  where $\eta<1-\beta.$
  Moreover if (A4) holds for all $Q\leq q\leq 2Q$, then for $X\geq Q^{2+\varepsilon}$ with $\log X\asymp \log Q$, we have
  \begin{equation}\label{eqn: lambdaF2mQ}
	\sum_{Q\leq q\leq 2Q}\sideset{}{}{\sum}_{{\chi\in\mathcal{F}_q}}\Big|\sum_{n\leq X}\lambda_F(n)\chi(n)\Big|^2\ll_F \frac{X^2}{\log^\eta X},
  \end{equation}
  where $\eta<1-\beta.$
\end{proposition}
To prove the above proposition, our main tool is the following theorem, which we will quickly establish in Section \S \ref{sec: tool thm}.
\begin{theorem}\label{thm: 2momentestimate}
Let $X, Q\geq 2$ and $Q\leq q\leq 2Q$, $\mathcal{E}_q$ be a subset of
$$\{\chi: \chi \text{ is a primitive Dirichlet character modulo } q\}.$$
Let $0\leq \beta<1$ and
let $f:\mathbb{N}\rightarrow \mathbb{C}$ be a multiplicative function which satisfying the following conditions:\\
	\begin{itemize}
		\item[C1]
		(Divisor bound) For any integer $n\geq 1$,
		$f(n)\ll \tau_d(n)$ for some fixed constant $d\geq 2$.\\
		\item[C2]
		(Average bound) For $X^{1-o(1)}\leq Y\leq X$, we have
		\begin{equation*}
			\sum_{X\leq n\leq X+Y}|f(n)|^2\ll Y.
		\end{equation*}
		\item[C3]
		(Sieve bound) For any $2\leq P_1\leq P_2\leq X$, we have
		\begin{equation*}
			\sum_{n\leq X\atop (n, \prod_{P_1\leq p\leq P_2}p)=1}|f(n)|^2\ll X\prod_{P_1\leq p\leq P_2}\Big(1-\frac{1}{p}\Big).
		\end{equation*}
		\item[C4]
		(Cancellation over primes)
        For $\exp(\log^{\beta+\varepsilon} X) \leq Y\leq \exp(\frac{\log X}{\log\log X}) $ and $\frac{Y}{\log X}\leq h\leq Y$, we have,
		\begin{equation*}
			\sup_{\chi\in \mathcal{E}_q}\Big|\sum_{Y\leq p\leq Y+h}f(p)\chi(p)\Big|\ll \frac{Y}{(\log Y)^{1+\gamma}},
		\end{equation*}
		for some $\gamma>0$.
	\end{itemize}
	Then we have
	\begin{equation}\label{eqn: q_2m}
		\sideset{}{}{\sum}_{\chi\in \mathcal{E}_q}\Big|\sum_{n\leq X}f(n)\chi(n)\Big|^2\ll qX^{1+\varepsilon}+\frac{X^2}{\log^\eta X},
	\end{equation}
	for any $\eta<\min(\frac{2}{3}\gamma,1-\beta).$
	Moreover, if C4 holds for all $Q\leq q\leq 2Q$, then we have
	\begin{equation}\label{eqn: Q_2m}
		\sum_{Q\leq q\leq 2Q}\sideset{}{}{\sum}_{\chi\in \mathcal{E}_q}\Big|\sum_{n\leq X}f(n)\chi(n)\Big|^2\ll Q^2X^{1+\varepsilon}+\frac{X^2}{\log^\eta X},
	\end{equation}
	for any $\eta<\min(\frac{2}{3}\gamma,1-\beta).$
\end{theorem}
\begin{remark}\label{rmk:p loglog bound}
	Conditions \emph{C2} and \emph{C3} can be replaced by the following bound
	\begin{equation}\label{eqn: p_loglog bound}
		\sum_{p\leq X}\frac{|f(p)|^2}{p}= \log\log X+O_f(1),
	\end{equation}
for all $X\geq 2$.
	Indeed, providing \eqref{eqn: p_loglog bound}, by Shiu's theorem \cite{Shiu1980}, we have
	\begin{equation}
		\sum_{X\leq n\leq X+Y}|f(n)|^2\ll Y\prod_{p\leq X}\Big(1+\frac{|f(p)|^2-1}{p}\Big)\ll Y.
	\end{equation}
	And by Lemma \ref{lem: achieve C3}, we can achieve \emph{C3}.
\end{remark}
\begin{remark}\label{rmk: zerofreeregion} Condition \emph{C4}  is somewhat difficult if $q$ is large. It is a variant of prime number theorem for $f\times\chi$. Define the twisted $L-$function associated $f$ by
	\begin{equation*}
		L(s, f\times \chi):=\sum_{n\geq 1}\frac{f(n) \chi(n)}{n^s}, \quad \text{ for } \Re(s)>1,
	\end{equation*}
	where $\chi$ is a primitive Dirichlet character modulo $q$.
	Assume that each $L(s, f\times \chi)$ has a nice analytic property (e.g. $f=F$ in Section \S \ref{subsec:mresult}), then \emph{C4} can be induced from the following zero free region:
	\begin{equation}
		\prod_{\chi\in \mathcal{E}_q}L(s, f\times \chi)\neq 0, \quad \text{ for } \Re(s)>1-\frac{c_f}{\log^{\beta}(q(3+|\Im(s)|))},
	\end{equation}
	where $0\leq\beta<1$ and $c_f$ is a positive constant depending only on $f$. \revisedtwo{This corresponds to the assumption (\emph{A4})} in Section \S \ref{subsec:mresult}.
The details of the argument can be find in Section \S \ref{sec:Zfr&pnt}.
\end{remark}
Theorem \ref{thm: 2momentestimate} can be compared with the large sieve inequalities
\begin{equation}\label{eqn: q-ls}
	\sideset{}{^*}{\sum}_{\chi (\mod q)}\Big|\sum_{n\leq X}a(n)\chi(n)\Big|^2\ll (q+X)\sum_{n\leq X}|a_n|^2,
\end{equation}
and
\begin{equation}\label{eqn: Q-ls}
	\sum_{Q\leq q\leq 2Q}\sideset{}{^*}{\sum}_{\chi (\mod q)}\Big|\sum_{n\leq X}a(n)\chi(n)\Big|^2\ll (Q^2+X)\sum_{n\leq X}|a_n|^2.
\end{equation}
These are classic $L^2$-$L^2$ type estimates and the norm $q+X$ and $Q^2+X$ are essentially best respectively.
Here $q$ and $Q^2$ come from the diagonal contributions which are unable to improve.
In our case, providing $a_n$ with multiplicative structure and some arithmetical conditions,
we can obtain a slight saving on the size of length which comes from the off diagonal contribution.\par
Such a second moment estimate is not surprising for the (non-pretentious) multiplicative functions. Similar to \eqref{eqn: q-ls} and \eqref{eqn: Q-ls}, we have the following well-known estimate of integral form
\[
\int_{0}^{T}\Big|\sum_{X< n\leq 2X}\frac{a_n}{n^{1+it}}\Big|^2\dd t\ll(T+X)\sum_{X< n\leq 2X}\frac{|a_n|^2}{n^2}.
\]
\revised{Let $\delta>0$ be given.} In \cite{M-R15}, using the Vinogradov--Korobov type zero-free region for the Riemann $\zeta$-function, Matom\"aki and Radziwi\l\l\, showed that
\begin{equation}
  \int_{0}^{T}\Big|\sum_{X< n\leq 2X}\frac{\lambda(n)}{n^{1+it}}\Big|^2\dd t\ll \frac{1}{(\log X)^{\frac{1}{3}-o(1)}}\Big(\frac{T}{X}+1\Big)+\frac{T}{X^{1-\revised{\delta/2}}},
\end{equation}
where $\lambda(n)$ is the Liouville function.
For the general $1$-bounded multiplicative function $f$, based on Hal\'asz's theorem, they showed that \revised{(}\revised{see also} \cite{M-R16}\revised{)},
\begin{equation}
  \int_{(\log X)^{\frac{1}{15}}}^{T}\Big|\sum_{X< n\leq 2X}\frac{f(n)}{n^{1+it}}\Big|^2\dd t\ll \frac{1}{(\log X)^{\frac{1}{48}}}\Big(\frac{T}{X}+1\Big)+\frac{T}{X^{1-o(1)}}.
\end{equation}
\revised{Actually, the above two estimates are nontrival for $X\geq T^{1+o(1)}$, since in this case we get the log-saving from the first term.}
\eqref{eqn: q_2m} and \eqref{eqn: Q_2m} can be viewed as the variants of the above two estimates by averaging the multiplicative characters $\chi$ instead of $n^{it}$. And our proof of Theorem \ref{thm: 2momentestimate} obeys the former case since that we can not obtain the non-trivial estimate for the mean value of $f\chi$ unconditionally.

\medskip
\textbf{Notation.}
Throughout the paper, $\varepsilon$ is an arbitrarily small positive number and $A$ is an arbitrarily large positive number,
all of them may be different at each occurrence.
As usual, $p$ stands for a prime number.

We use the standard Landau and Vinogradov asymptotic notations \( O(\cdot) \), \( o(\cdot) \), \(\ll\), and \(\gg\). Specifically, we express \( X \ll Y \), \( X = O(Y) \), or \( Y \gg X \) when there exists a constant \( C \) such that \( |X| \leq CY \). If the constant $C=C_s$ depends on some object $s$, we write $X=O_s(Y)$. And the notation \( X \asymp Y \) is used when both \( X \ll Y \) and \( Y \ll X \) hold. As \( N \to \infty \), \( X = o(Y) \) indicates that \( |X| \leq c(N)Y \) for some function \( c(N) \) that tends to zero.

For a proposition \( P \), we write \( 1_P(x) \) for its indicator function, i.e. a function that equals to 1  if \( P \) is true and 0 if \( P \) is false.

\section{ Some bounds \revised{on} multiplicative functions}
The first lemma is given by Tenenbaum in \cite[\revised{Chapter III, Theorem 3.5}]{Tenenbaum15},  which provides an upper bound for the average of a class of general multiplicative functions in terms of their logarithmic average. Here we use the \revised{qualitative} version.

\begin{lemma}\label{lem: Hall-T}
	Let \( f \) be a non-negative multiplicative function such that
	\[
	\sum_{p \leq y} f(p) \log p \ll y \quad (y \geq 0),
	\]
	and\[
	\sum_{p} \sum_{\nu \geq 2} \frac{f(p^\nu)}{p^\nu} \log p^\nu \ll 1.
	\]
	Then, for \( X > 1 \), we have
	\[
	\sum_{n \leq X} f(n) \ll \frac{X}{\log X} \sum_{n \leq X} \frac{f(n)}{n}.
	\]
\end{lemma}
Notice that we have, on the right,
\[
\sum_{n \leq X} \frac{f(n)}{n} \leq \prod_{p \leq X} \left( 1 + \frac{f(p)}{p} + \frac{f(p^2)}{p^2} + \cdots \right),
\]
the form in which the result is often applied.  As an application of Lemma \ref{lem: Hall-T}, let's see how to achieve  \emph{C3} by using the bound \eqref{eqn: p_loglog bound}.
\begin{lemma}\label{lem: achieve C3}
	Assume that $f(n)$ is a \(\tau_d\)-bounded multiplicative function and the bound \eqref{eqn: p_loglog bound} holds. We have
\begin{equation*}
		\sum_{n\leq X\atop (n, \prod_{P_1\leq p\leq P_2}p)=1}|f(n)|^2\ll X\prod_{P_1\leq p\leq P_2}\Big(1-\frac{1}{p}\Big)
	\end{equation*}
uniformly for $2\leq P_1\leq P_2\leq X$ .
\end{lemma}
\begin{proof}
	Clearly, $|f(n)|^2$ is multiplicative and non-negative. Suppose that for any integer $n\geq 1$,
	$f(n)\ll \tau_d(n)$ for some fixed constant $d\geq 2$.  We verify the conditions in Lemma \ref{lem: Hall-T}.
	
	Using the bound \eqref{eqn: p_loglog bound} and summation by parts, we obtain
	\begin{equation*}
		\sum_{p \leqslant X}\left|f(p)\right|^{2} \log p \ll X .
	\end{equation*}
	And we have
	\[\sum_{p} \sum_{\nu \geq 2} \frac{f(p^\nu)}{p^\nu} \log p^\nu \ll 1,\]
	since $f(n)$ is divisor bounded. Hence by Lemma \ref{lem: Hall-T}, we get
	\revisedtwo{	\begin{align}
			\sum_{n\leq X\atop (n, \prod_{P_1\leq p\leq P_2}p)=1}|f(n)|^2
			&\ll  \frac{X}{\log X}\sum_{n\leq X\atop (n, \prod_{P_1\leq p\leq P_2}p)=1}\frac{|f(n)|^2}{n} \nonumber \\ 
			&\ll \frac{X}{\log X}\prod_{\substack{p \leqslant X \\ p \notin [P_1,P_2]}}\left(1+\frac{|f(p)|^2}{p}+\sum_{\nu \geqslant 2} \frac{f\left(p^{\nu}\right) }{p^{\nu}}\right) \nonumber\\
			&\ll \frac{X}{\log X} \exp \bigg(\sum_{\substack{p \leqslant X \\ p \notin [P_1,P_2]}}\frac{|f(p)|^2}{p}+O(1)\bigg).\label{eqn: exp bpund}
		\end{align}}
	By the bound \eqref{eqn: p_loglog bound}, we derive
	\begin{equation}\label{eqn: 3 loglog}
		\begin{aligned}
			\sum_{\substack{p \leqslant X \\
					p \notin [P_1,P_2]}} \frac{|f(p)|^2}{p} & =\sum_{p \leqslant X} \frac{|f(p)|^2}{p}-\sum_{p \leqslant P_2} \frac{|f(p)|^2}{p}+\sum_{p \leqslant P_1} \frac{|f(p)|^2}{p} \\
			& =\log \log X-\log \log P_2+\log \log P_1+O(1).
		\end{aligned}
	\end{equation}
	Then the claim follows by inserting \eqref{eqn: 3 loglog} into \eqref{eqn: exp bpund} and Mertens' theorem.
\end{proof}

\section{Zero free region and prime number theorem}\label{sec:Zfr&pnt}
In this section, we consider the relationship between the zero free region with the prime number theorem on $\lambda_F \chi$.
This corresponds {to} what we have claimed in Remark \ref{rmk: zerofreeregion}.
The key argument is the explicit formula for summing $\lambda_F\chi$ over primes; see \eqref{eqn: explicitformula}.\par

\revisedtwo{Recall the definition \eqref{eqn:L F twist chi}
and note that $L(s,F\time \chi)$ is absolutely convergent} in $\Re(s)>1.$ For the sake of convenience, we  write
\begin{equation}\label{eqn: coeef a_F}
	a_F(p^k)=\sum_{j=1}^{d}\alpha_{j,F}(p)^k,
\end{equation}
and set $a_F(n)=0$ if $n$ is not a prime power\revisedtwo{. We have}
\begin{equation}
	\log L(s,F\times \chi)=\sum_{n\geq 2}\frac{\Lambda(n)a_F(n)\chi(n)}{n^s \log n},
\end{equation}
and
\begin{equation}\label{eqn: log deri}
	-\frac{L'}{L}(s,F\times \chi)=\sum_{n\geq 1}\frac{\Lambda(n)a_F(n)\chi(n)}{n^s}.
\end{equation}
Note that $ a_F(p)=\lambda_F(p)$ for $p$ prime.

By the absolute convergence of Euler product for $\Lambda(s, F\times\chi)$ and \eqref{eqn: L_infty} with non-vanishing of $\Gamma$-functions, we see that there is no zero in the half plane $\Re(s)>1$. By the functional equation \eqref{eqn: FE},  the same is true in the half plane $\Re(s)<0$. So the zeros of $\Lambda(s, F\times\chi)$ must lie in the in the critical strip $0\leq \Re(s)\leq 1$, and we call them the non-trivial zeros of $L(s, F\times\chi)$. We also make the assumption that neither 0 nor 1 can be a non-trivial zero of $L(s,F\times\chi)$.

By hypothesis, the complete $L$-function $\Lambda(s, F\times\chi)$ is an entire function of order 1, and thus there exist constants $a = a(F,\chi)$ and $b = b(F,\chi)$ such that it has a Hadamard product representation
\begin{equation}
	L(s, F\times \chi)L_{\infty}(s, F\times \chi)=e^{a+bs}\prod_{\rho}(1-\frac{s}{\rho})e^{s/\rho},
\end{equation}
where $\rho$ runs through the non-trivial zeros mentioned above. By taking the logarithmic derivative of both sides, we have
\begin{equation}\label{eqn: log deri fac}
	\frac{L'}{L}(s,F\times\chi)+\frac{L'_{\infty}}{L_{\infty}}(s,F\times\chi)=b+\sum_{\rho}\left(\frac{1}{s-\rho}+\frac{1}{\rho}\right).
\end{equation}
Using  the functional equation \eqref{eqn: FE} and the fact that $ L(s, F\times \chi)L_{\infty}(s, F\times \chi)$ is an entire function of order 1, one can prove that $\Re(b)$ equals the absolutely convergent sum $-\sum_{\rho}\Re(\rho^{-1})$.
It follows that
\begin{equation}\label{eqn: real part}
	\sum_{\rho}\Re\left(\frac{1}{s-\rho}\right)=\Re \left(\frac{L'}{L}(s,F\times\chi)+\frac{L'_{\infty}}{L_{\infty}}(s,F\times\chi) \right).
\end{equation}
\begin{lemma}\label{lem: number of zeros}
	Let $T\geq 0$. \revisedtwo{Let $\rho$ denote the non-trivial zeros} of $L(s, F\times\chi)$. Then
	\begin{equation}\label{eqn: non-trivial zeros recip}
		\sum_{\rho}\frac{1}{1+(T-\Im\rho)^2}\ll_F \log \big(q(T + 2)\big).
	\end{equation}
	Moreover, the number  $m(T, F \times \chi)$  of non-trivial zeros  $\rho$  such that  $|\Im\rho-T| \leqslant 1$  is $O_F(\log \big(q(T + 2)\big) )$.
\end{lemma}
\begin{proof}
	Taking the logarithmic derivative on both sides of equation \eqref{eqn: L_infty}, we get
	\begin{equation}
		\frac{L'_{\infty}}{L_{\infty}}(s,F\times\chi)=-\frac{d}{2}\log\pi+\frac{1}{2}\log N+\frac{1}{2}\sum_{j=1}^d \frac{\Gamma^{\prime}(\frac{s+\mu_j}{2} )}{\Gamma(\frac{s+\mu_j}{2} )}.
	\end{equation}
	By definition
	\[
	\frac{1}{\Gamma(s)}=s e^{\gamma s} \prod_{n=1}^{\infty}\left(1+\frac{s}{n}\right) e^{-\frac{s}{n}} .
	\]
	Thus
	\[
	-\frac{1}{2} \frac{\Gamma^{\prime}\left(\frac{s+\mu_j}{2}\right)}{\Gamma\left(\frac{s+\mu_j}{2}\right)}=\frac{1}{s+\mu_j}+\frac{\gamma}{2}+\sum_{n\geq 1}\left(\frac{1}{2 n+s+\mu_j}-\frac{1}{2 n}\right) .
	\]
	Apply   \eqref{eqn: log deri fac}, we get
	\begin{equation}\label{eqn: L'/L factori}
		\begin{aligned}
			\frac{L'}{L}(s,F\times\chi)=&b+\frac{d}{2}\log\pi-\frac{1}{2}\log N+\frac{d}{2}\gamma+\sum_{j=1}^d \frac{1}{s+\mu_j}  \\&+\sum_{\rho}\left(\frac{1}{s-\rho}+\frac{1}{\rho}\right)+\sum_{j=1}^d \sum_{n\geq 1}\left(\frac{1}{2 n+s+\mu_j}-\frac{1}{2 n}\right).
		\end{aligned}
	\end{equation}
	Note that for $s=\sigma+it$, $-1\leq\sigma\leq 2$,
	\begin{equation}
		\left|\sum_{j=1}^d \sum_{n\geq1}\left(\frac{1}{2 n+s+\mu_j}-\frac{1}{ n}\right) \right|\ll d \sum_{n\leq |t|}\frac{1}{2n}+d\sum_{n>|t|}\frac{|s|}{4n^2}\ll_F d\log (|t|+2).
	\end{equation}
	Hence by \eqref{eqn: real part}, we have
	\begin{equation}
		\Re \frac{L'}{L}(s,F\times\chi)=\Re\sum_{j=1}^{d}\frac{1}{s+\mu_j}+\Re \sum_{\rho} \frac{1}{s-\rho}+O_F(d \log \big(q(T + 2)\big)).
	\end{equation}
	Let  $s=2+i T $. Observe that
	\begin{equation}
		\Re \sum_{\rho}\frac{1}{2+iT-\rho}=\sum_{\rho}\frac{2-\Re\rho}{(2-\Re\rho)^2+(T-\Im\rho)^2}\gg \sum_{\rho}\frac{1}{1+(T-\Im\rho)^2}.
	\end{equation}
	Thus, \eqref{eqn: non-trivial zeros recip} follows from the estimate
	\begin{equation}
		\begin{aligned}
				\Re \sum_{\rho}\frac{1}{2+iT-\rho}&\ll_F d \log \big(q(T + 2)\big)+\Big|\frac{L'}{L}(2+iT,F\times\chi) \Big|\\
				&\ll_F d\log \big(q(T + 2)\big).
		\end{aligned}
	\end{equation}
	This immediately implies that the number  $m(T, F \times \chi)$  of non-trivial zeros  $\rho$  such that  $|\Im\rho-T| \leqslant 1$  is
	\begin{equation}
		\sum_{|\Im\rho-T|\leq 1}1\leq 2\sum_{\rho}\frac{1}{1+(T-\Im\rho)^2}\ll_F d\log \big(q(T + 2)\big),
	\end{equation}
	where the first summation runs over all non-trivial zeros $\rho$ satisfying $|\Im\rho-T|\leq 1$.
\end{proof}
\begin{lemma}\label{lem: L'/L esti}
	Let $\rho$ denote the non-trivial zeros of $L(s, F\times\chi)$. For $s=\sigma+it$ with $-1\leq\sigma\leq 2$, we have
	\begin{equation}\label{eqn: L'/L esti}
		\frac{L^{\prime}}{L}(s,F\times\chi)=\sum_{1\leq j\leq d\atop|s+\mu_{j}|\leq 1} \frac{1}{s+\mu_{j}}+\sum_{|t-\Im\rho|\leq 1} \frac{1}{s-\rho} +O_F\left(  d \log \big(q(|t| + 2)\big)  \right),
	\end{equation}
	where the second summation runs over all non-trivial zeros $\rho$ satisfying $|\Im\rho-t|\leq 1$.
\end{lemma}
\begin{proof}
	By \eqref{eqn: L'/L factori}, we have
	\begin{equation}\label{eqn: L'/L diff}
		\begin{aligned}
			\frac{L^{\prime}}{L}(s,F\times\chi)&-\frac{L^{\prime}}{L}(2+it,F\times\chi)=\sum_{j=1}^d\frac{1}{s+\mu_j} \\&+\sum_{\rho}\left(\frac{1}{\sigma+it-\rho}-\frac{1}{2+it-\rho}\right)
			+O_F\left(d \log \big(q(|t| + 2)\big)  \right).
		\end{aligned}
	\end{equation}
	Note that \(\frac{L^{\prime}}{L}(2+it,F\times\chi)=O_F(1)\). Then for \(|t-\Im\rho|>1\),
	\[
	\left|\frac{1}{\sigma+it-\rho}-\frac{1}{2+it-\rho} \right|\leq \frac{3}{|t-\Im\rho|^2}\leq \frac{6}{1+|t-\Im\rho|^2}.
	\]
	This together with Lemma \ref{lem: number of zeros} show that the sum over $\rho$ with $|t-\Im\rho|\geq 1$ in \eqref{eqn: L'/L diff} is $O_F\left(d\log q(|t|+2)\right)$. Moreover, for $|t-\Im\rho|\leq 1$,
	\[
	\sum_{|t-\Im\rho|\leq 1}\left|\frac{1}{2+it-\rho}\right|\ll\sum_{|t-\Im\rho|\leq 1} 1\ll_F d \log \big(q(|t| + 2)\big)  .
	\]
	The desired result \eqref{eqn: L'/L esti} follows from \eqref{eqn: L'/L diff}.
\end{proof}
\begin{lemma}\label{lemma: explicitformula}
	\revisedtwo{Let $\rho$ denote the non-trivial zeros} of $L(s, F\times\chi)$. Assuming GRC on $F$, for $X, T, q\geq 2$, we have
	\begin{equation}\label{eqn: explicitformula}
		\sum_{p\leq X}\lambda_F(p)\chi(p)\log p=-\sum_{|\Im\rho|\leq T}\frac{X^\rho}{\rho}+O_F\Big(\frac{X}{T}(\log X\log(qT))^{2}+X^{\frac{1}{2}+\varepsilon}\Big),
	\end{equation}
	where $\sum_{|\Im\rho|\leq T}$ means the summand of non-trivial zeros $\rho$ satisfying $|\Im\rho|\leq T$ with multiplicity.
\end{lemma}
\begin{proof}
	We denote
	\begin{equation}
		\psi(F\times\chi,X)=\sum_{n\leq X}\Lambda(n)a_F(n)\chi(n).
	\end{equation}
	Then we can see that
	\begin{equation}\label{eqn: psi error}
		\begin{aligned}
			\psi(F\times\chi,X)&= \sum_{p\leq X}a_F(p)\chi(p)\log p +\sum_{{n\leq X\atop n=p^k,k\geq2}}\Lambda(n)a_F(n)\chi(n)\\
			&=\sum_{p\leq X}\lambda_F(p)\chi(p)\log p+O_F\left(X^{1/2+\varepsilon}\right),
		\end{aligned}
	\end{equation}
	since $|a_F(n)|\leq d$ on GRC.
	It suffices to show
	\[
	\psi(F\times\chi,X)=-\sum_{|\Im\rho|\leq T}\frac{X^\rho}{\rho}+O\Big(\frac{X}{T}(\log X\log(qT))^{2}\Big).
	\]
	We can assume that fractional parts of  $X$ is  $1/2$. Perron's formula says that, for any  $c>0$  and  $y>0 $, we have
	\begin{equation*}
		\frac{1}{2 \pi i} \int_{c-i T}^{c+i T} y^{s} \cdot \frac{d s}{s}=\left\{\begin{array}{ll}
			1, & \text { if } y>1 \\
			0, & \text { if } y<1
		\end{array}+O\left(\frac{y^{c}}{\max (1, T|\log y|)}\right) .\right.
	\end{equation*}
	Therefore, by using \eqref{eqn: log deri} and letting  $c=1+1 / \log X $, we have
	\begin{equation}\label{perron1}
		\psi(F\times\chi,X)=\frac{1}{2 \pi i} \int_{c-i T_1}^{c+i T_1}-\frac{L'}{L}(s,F\times \chi) \cdot \frac{X^s}{s} d s+O_F\left(\frac{X\log^2 X}{T_1} \right),
	\end{equation}
	where $T\leq T_1\leq T+1$ such that the distance of nearest non-trivial zero of  \(L(s,F\times\chi)\) to the line \(\Im s=T_1\) is \(\gg_F 1/\log qT\), and the distance of nearest \(\Im\mu_j\) to \(T_1\) is \(\gg_F 1/\log^2 qT\).
	
	Note that the complete $L$-function $\Lambda(s, F\times\chi)$ is an entire function of order $1$. Shift the integral along a rectangular contour to the line $-1/2-\xi\pm iT_1$, where \(0\leq \xi< 1/2\) is a small constant such that the distance of nearest \(\Re\mu_j\) to \(-1/2-\xi\) is \(\gg_F 1/\log qT\)\revisedtwo{. Picking up residues} at all non-trivial zeros $\rho$ we get
	\begin{equation}\label{eqn: psi error_1}
		\psi(F\times\chi,X)=-\sum_{|\Im\rho|\leq T}\frac{X^\rho}{\rho}+R+O_F\left(\frac{X\log^2 X}{T_1} \right),
	\end{equation}
	where
	\begin{equation}\label{eqn: R bound}
		\begin{aligned}
			R
			& =-\frac{1}{2\pi i}\left(\int_{-\frac{1}{2}-\xi-iT_1}^{-\frac{1}{2}-\xi+iT_1}+\int_{-\frac{1}{2}-\xi+iT_1}^{c+iT_1}+\int_{c-iT_1}^{-\frac{1}{2}-\xi-iT_1}\right)\frac{L'}{L}(s,F\times\chi) \frac{X^{s}}{s} ds \\
			& \ll X^{-\frac{1}{2}-\xi}\int_{-T_1}^{T_1}\Big| \frac{L'}{L}\left(-\frac{1}{2}-\xi+it,F\times\chi\right)\Big|\frac{dt}{1+|t|}\\&\quad +\frac{X}{T_1}\int_{-\frac{1}{2}-\xi}^{c}\Big| \frac{L'}{L}\left(\sigma\pm iT_1,F\times\chi\right)\Big| d\sigma.
		\end{aligned}
	\end{equation}
	By Lemma \ref{lem: number of zeros} and Lemma \ref{lem: L'/L esti}, 
	for $s=\sigma\pm iT_1$ with \(-1/2-\xi\leq\sigma\leq c\) and \(|T_1-\Im\mu_j|\gg_F 1/ (\log^2 qT) \) for each \(1\leq j\leq d\), we have
	\[
	\begin{aligned}
		\Big| \frac{L'}{L}\left(\sigma\pm iT_1,F\times\chi\right)\Big|&= \sum_{1\leq j\leq d\atop|s+\mu_{j}|\leq 1} \frac{1}{\left|(\sigma+\Re\mu_{j})+ i(\Im\mu_j\pm T_1)\right|}\\
		&+ \sum_{|T_1-\Im\rho|<1} \frac{1}{|\sigma-\Re\rho+i(T_1-\Im\rho)|} +O_F\left(d\log qT\right)\\
		&\ll_F (\log^2 qT)+ \sum_{|T_1-\Im\rho|<1}\frac{1}{|T_1-\Im\rho|}\ll_F (\log qT)^2.
	\end{aligned}
	\]
	Here we have used the fact that \(\left|T_1-\Im\rho \right|\gg_F 1/\log qT\).  For $|t|\leq T_1$, similarly we have
	\[
	\begin{split}
		\Big| \frac{L'}{L}\left(-\frac{1}{2}-\xi+it,F\times\chi\right)\Big|
		&= \sum_{1\leq j\leq d\atop|s+\mu_{j}|\leq 1} \frac{1}{\left|(-\frac{1}{2}-\xi+\Re\mu_{j})+ i(\Im\mu_j\pm t)\right|}\\
		&+ \sum_{|t-\Im\rho|<1} \frac{1}{|-\frac{1}{2}-\xi-\Re\rho+i(t-\Im\rho)|} \\
		&+O_F\left(d\log q(|t|+2)\right)\\
		&\ll_F \log qT.
	\end{split}
	\]
	Hence by \eqref{eqn: R bound}
	\begin{equation}\label{eqn: R final bound}
		R\ll_F \left(X^{-1/2}+\frac{X}{T}\right)\log^2 qT.
	\end{equation}
	The desired result \eqref{eqn: explicitformula} follows by combining \eqref{eqn: R final bound}, \eqref{eqn: psi error_1} and \eqref{eqn: psi error}.
\end{proof}
We conclude this section by explaining how to derive condition (\emph{C4}), \revisedtwo{actually} a stronger form, from a good zero free region for $L(s,F\times\chi)$, as discussed  in Remark \ref{rmk: zerofreeregion}.
\begin{lemma}\label{lemma: bound p-sum}
	Let $Y, q\geq 2$.
    Assuming GRC and the following zero free region:
	\begin{equation}\label{eqn: ZFR}
		 L(\sigma+it, F\times \chi)\neq 0, \quad \text{ for } \sigma>1-\frac{c_F}{\log^{\beta}(qT)},\quad |t|\leq T,
	\end{equation}
    where $0<\beta<1$, $\log^{100}q \ll T \ll q^{o(1)}$, and $c_F$ is a positive constant depending only on $F$. Then for any $Y^{\frac{1}{2024}}\leq h\leq Y$, we have,
    uniformly for $\log^{\beta+\varepsilon}q\ll \log Y\ll \log^A q$, the bound
	\begin{equation*}
		\Big|\sum_{Y\leq p\leq Y+h}\lambda_F(p)\chi(p)\Big|\ll_{F} \frac{Y}{\log^{A}Y}
	\end{equation*}
    holds for arbitrarily large positive number $A$, and the $\ll$ implied constant may depend on $F$ but it is independent of $\chi$.
    Moreover if \eqref{eqn: ZFR} holds for all $L(s, F\times \chi)$ with $\chi \in \mathcal{E}_q$, where $\mathcal{E}_q$ is subset of all primitive Dirichlet characters of conductor $q$,
    under the same conditions, we have
    \begin{equation*}
		\sup_{\chi\in \mathcal{E}_q}\Big|\sum_{Y\leq p\leq Y+h}\lambda_F(p)\chi(p)\Big|\ll_{F} \frac{Y}{\log^{A}Y}.
	\end{equation*}
\end{lemma}
\begin{proof}
From GRC and the Brun--Titchmarsh theorem, we have
\begin{equation*}
  \Big|\sum_{Y\leq p\leq Y+h}\lambda_F(p)\chi(p)\Big|\ll \sum_{Y\leq p\leq Y+h}1\ll\frac{h}{\log Y}.
\end{equation*}
Thus it suffices to consider $\frac{Y}{\log^{A}Y}\ll h\leq Y.$
Let $T\geq 3$ be a parameter which will be chosen later.
Assuming the zero free region for $L(s, F\times \chi)$ \eqref{eqn: ZFR},  by Lemma \ref{lemma: explicitformula}, we have
\begin{equation}\label{eqn: f chi weight log}
		\begin{aligned}
			\Big| \sum_{Y\leq p\leq Y+h}&\lambda_F(p)\chi(p)\log p\Big|\\
            &\leq \Big|\sum_{|\Im \rho|\leq T}\frac{(Y+h)^{\rho}-Y^{\rho}}{\rho}\Big|+O_F\Big(\frac{Y}{T}(\log Y\log(qT))^{2}+Y^{\frac{1}{2}+\varepsilon}\Big)\\
			&\leq \int_{Y}^{Y+h}\sum_{|\Im \rho|\leq T}y^{\Re\rho-1}\dd y+O_F\left(\frac{Y}{T}(\log Y\log(qT))^{2}+Y^{\frac{1}{2}+\varepsilon}\right)\\
            &\ll_F hT\log qT\exp\Big(-c_F\frac{\log Y}{\log^{\beta}qT}\Big)+\frac{Y}{T}(\log Y\log(qT))^{2}+Y^{\frac{1}{2}+\varepsilon}.
		\end{aligned}
	\end{equation}
Taking \(T=\log^A Y\), for $\log^{\beta+\varepsilon}q\ll \log Y\ll \log^A q$, the error term of right side of above is acceptable. And \eqref{eqn: f chi weight log} becomes
	\begin{equation*}
		\begin{aligned}
			\sum_{Y\leq p\leq Y+h}\lambda_F(p)\chi(p)\log p&\ll_{F} h(\log^A Y)\log qT \exp\left(-c_F\frac{\log Y}{\log^\beta  qT}\right)+\frac{Y}{\log^A Y}\\
			&\ll_{F} \frac{Y}{\log^A Y}
		\end{aligned}
	\end{equation*}
	for arbitrarily large positive number $A$.
    Applying the partial summation formula, we have
    \begin{multline*}
    \sum_{Y\leq p\leq Y+h}\lambda_F(p)\chi(p)\ll \frac{1}{\log(Y+h)}\Big|\sum_{Y\leq p\leq Y+h}\lambda_F(p)\chi(p)\log p\Big|\\
    +\int_{0}^{h}\Big|\sum_{Y\leq p\leq Y+y}\lambda_F(p)\chi(p)\log p\Big|\frac{\dd y}{Y+y}\ll \frac{Y}{\log^A Y}.
    \end{multline*}
     This completes the proof.
\end{proof}

\section{Bounds for second moment estimates }\label{sec: tool thm}
\subsection{Proof of Theorem \ref{thm: 2momentestimate}}
To prove Theorem \ref{thm: 2momentestimate}, we shall suppose that $q\sim Q\leq X$. Otherwise, it is trivially followed by using
\eqref{eqn: q-ls} and \eqref{eqn: Q-ls}.
The core of the proof of Theorem \ref{thm: 2momentestimate} is Ramar\'e's identity:
\begin{equation}\label{eqn: RamareId}
	\sum_{mp=n\atop P_1\leq p\leq P_2}\frac{1}{\omega_{[P_1, P_2]}(m)+1_{p\nmid m}}=
	\begin{cases}
		& 1, \quad \text{if} \quad (n, \prod_{P_1\leq p\leq P_2}p)\neq 1,\\
		& 0, \quad  \text{otherwise},\\
	\end{cases}
\end{equation}
for any fixed $n\geq 1$ and $2\leq P_1\leq P_2$ where  $\omega_{[P_1, P_2]}(m)$  is the number of distinct prime divisors of  $m$  on  $[P_1, P_2]$. This identity follows directly since  the number of representations  $n=p m$  with  $P_1\leq p \leq P_2$  is  $\omega_{[P_1, P_2]}(m)$. Using \eqref{eqn: RamareId}, we have
\begin{equation}\label{eqn: Dirichlet_f}
	\sum_{n\leq X}f(n)\chi(n)
	=\sum_{P_1\leq p\leq P_2}\sum_{m\leq X/p}\frac{f(pm)\chi(pm)}{\omega_{[P_1, P_2]}(m)+1_{p\nmid m}}+\sum_{n\leq X\atop (n, \prod_{P_1\leq p\leq P_2}p)= 1}f(n)\chi(n).
\end{equation}
By the multiplicity of $f\chi$, we rewrite the first sum as
\begin{multline}
	\sum_{P_1\leq p\leq P_2}f(p)\chi(p)\sum_{m\leq X/p}\frac{f(m)\chi(m)}{\omega_{[P_1, P_2]}(m)+1}
	\\+\sum_{P_1\leq p\leq P_2}\sum_{m\leq X/p\atop p\mid m}\Big(\frac{f(pm)\chi(pm)}{\omega_{[P_1, P_2]}(m)}-\frac{f(p)f(m)\chi(p)\chi(m)}{\omega_{[P_1, P_2]}(m)+1}\Big).
\end{multline}
Splitting the first sum above into small ranges, we  get that
\begin{equation}\label{eqn: decompostion_f}
	\sum_{\lfloor H\log P_1\rfloor\leq j\leq H\log P_2}\sum_{e^{j/H}< p\leq e^{(j+1)/H}\atop P_1\leq p\leq P_2}f(p)\chi(p)\sum_{m\leq Xe^{-j/H}\atop mp\leq X}\frac{f(m)\chi(m)}{\omega_{[P_1, P_2]}(m)+1},
\end{equation}
where
$\exp(\log^{\beta+\varepsilon} X)\leq P_1\leq P_2\leq \exp{\Big(\frac{\log X}{\log\log X}\Big)}$ and $H\gg \log\log X$ will be chosen later.
We remove the condition $mp\leq X$ over-counting at most by the integers $mp$
in the range $(X, Xe^{1/H}]$. Therefore we can, for some sequences $d_n$ bounded by
\begin{equation}
	|d_n|\leq \sum_{mp=n\atop  P_1\leq p\leq P_2}\Big|\frac{f(p)f(m)}{\omega_{[P_1,P_2]}(m)+1}\Big|,
\end{equation}
rewrite \eqref{eqn: decompostion_f} as
\begin{align*}
	\sum_{\lfloor H\log P_1\rfloor\leq j\leq H\log P_2}\sum_{e^{j/H}< p\leq e^{(j+1)/H}\atop P_1\leq p\leq P_2}f(p)\chi(p)\sum_{m\leq Xe^{-j/H}}\frac{f(m)\chi(m)}{\omega_{[P_1, P_2]}(m)+1}\\
	 +\sum_{X< n\leq Xe^{1/H}}d_n\chi(n).\qquad
\end{align*}
Thus, \eqref{eqn: Dirichlet_f}  is equal to
\begin{multline}\label{eqn: decomp}
	\sum_{\lfloor H\log P_1\rfloor\leq j\leq H\log P_2}\sum_{e^{j/H}<p\leq e^{(j+1)/H}\atop P_1\leq p\leq P_2}f(p)\chi(p)\sum_{m\leq Xe^{-j/H}}\frac{f(m)\chi(m)}{\omega_{[P_1,P_2]}(m)+1}\\
	+\sum_{P_1\leq p\leq P_2}\sum_{m\leq X/p\atop p\mid m}\Big(\frac{f(pm)\chi(pm)}{\omega_{[P_1, P_2]}(m)}-\frac{f(p)f(m)\chi(p)\chi(m)}{\omega_{[P_1, P_2]}(m)+1}\Big)\\+\sum_{X< n\leq Xe^{1/H}}d_n\chi(n)
	+
	\sum_{n\leq X\atop (n, \prod_{P_1\leq p\leq P_2}p)= 1}f(n)\chi(n).\qquad\quad
\end{multline}
Squaring this and taking summation over  $\chi \in\mathcal{E}_q$, then for some $j\in [\lfloor H\log P_1\rfloor, H\log P_2]$,
we have
\begin{multline}\label{eqn: 2mebound}
	\sideset{}{}{\sum}_{\chi\in\mathcal{E}_q}\Big|\sum_{n\leq X}f(n)\chi(n)\Big |^2
	\ll
	\Big(H\log \frac{P_2}{P_1}\Big)^2\sup_{\chi\in\mathcal{E}_q}\Big|\sum_{e^{j/H}< p\leq e^{(j+1)/H)}\atop P_1\leq p\leq P_2}f(p)\chi(p)\Big|^2
	\\ \times\sideset{}{^*}{\sum}_{\chi (\mod q)}\Big|\sum_{m\leq Xe^{-j/H}}\frac{f(m)\chi(m)}{\omega_{[P_1,P_2]}(m)+1}\Big|^2
	+ \sideset{}{^*}{\sum}_{\chi (\mod q)}\Big|\sum_{P_1\leq p\leq P_2}\sum_{m\leq X/p\atop p\mid m}\frac{f(pm)\chi(pm)}{\omega_{[P_1, P_2]}(m)}\Big|^2
	\\+\sideset{}{^*}{\sum}_{\chi (\mod q)}\Big|\sum_{P_1\leq p\leq P_2}\sum_{m\leq X/p\atop p\mid m}\frac{f(p)f(m)\chi(p)\chi(m)}{\omega_{[P_1, P_2]}(m)+1}\Big|^2+\sideset{}{^*}{\sum}_{\chi (\mod q)}\Big|\sum_{X< n\leq Xe^{1/H}}d_n\chi(n)\Big|^2
	\\+\sideset{}{^*}{\sum}_{\chi (\mod q)}\Big|
	\sum_{n\leq X\atop (n, \prod_{P_1\leq p\leq P_2}p)= 1}f(n)\chi(n)\Big|^2.\qquad\qquad
\end{multline}
For the first term of the right-hand side of \eqref{eqn: 2mebound}, applying \emph{C4} which gives
\[
\sup_{\chi\in\mathcal{E}_q}\Big|\sum_{e^{j/H}< p\leq e^{(j+1)/H}\atop P_1\leq p\leq P_2}f(p)\chi(p)\Big|
\ll e^{(j+1)/H}(\frac{H}{j+1})^{1+\gamma}
\]
and \eqref{eqn: q-ls}, it is bounded by
\begin{equation*}
	\Big(H\log \frac{P_2}{P_1}\Big)^2e^{2(j+1)/H}(\frac{H}{j+1})^{2+2\gamma}(q+Xe^{-j/H})\sum_{m\leq Xe^{-j/H}}|f(m)|^2.
\end{equation*}
Using \emph{C2}, this is bounded by
\begin{align*}
	\Big(H\log \frac{P_2}{P_1}\Big)^2e^{2(j+1)/H}(\frac{H}{j+1})^{2+2\gamma}(q+Xe^{-j/H}) Xe^{-j/H}\\
	\ll qX^{1+\varepsilon}+X^2\frac{H^{2}}{(\log P_2)^{2\gamma}}.
\end{align*}
The second sum and the third sum, by \eqref{eqn: q-ls} and \emph{C1}, contribute at most
\begin{multline*}
	(q+X)\sum_{mp^2\leq X\atop P_1\leq p\leq P_2}\tau_d^2(mp^2)\ll (q+X)\sum_{P_1\leq p\leq P_2}\sum_{m\leq X/p^2}\tau_d^2(m)\\
	\ll (qX+X^2)\frac{\log^{d^2-1} X}{P_1}\ll (qX+X^2)(\log X)^{-2024}.
\end{multline*}
The fourth sum, by \eqref{eqn: q-ls} and \emph{C2}, contributes at most
\begin{equation*}
\begin{split}
	(q+X)&\sum_{X< n\leq Xe^{1/H}}\Big(\sum_{mp=n\atop  P_1\leq p\leq P_2}\Big|\frac{f(p)f(m)}{\omega_{[P_1,P_2]}(m)+1}\Big|\Big)^2\\
    &\ll (q+X)\sum_{X< n\leq Xe^{1/H}}\omega_{[P_1,P_2]}(n)^2\sum_{mp=n\atop  P_1\leq p\leq P_2}\Big|\frac{f(p)f(m)}{\omega_{[P_1,P_2]}(m)+1}\Big|^2\\
	&\ll (q+X)\sum_{X< n\leq Xe^{1/H}}\sum_{n=pm\atop P_1\leq p\leq P_2}|f(m)|^2\\
	&\ll (q+X)\sum_{P_1\leq p\leq P_2}\sum_{X/p<m\leq Xe^{1/H}/p}|f(m)|^2\\
	&\ll \frac{qX+ X^2}{H}\log \log P_2.
\end{split}
\end{equation*}
The last sum, by \eqref{eqn: q-ls} and \emph{C3}, contributes at most,
\begin{equation}
	(q+X)\sum_{n\leq X\atop (n, \prod_{P_1\leq p\leq P_2}p)= 1}|{f(n)}|^2\ll (qX+ X^2)\frac{\log P_1}{\log P_2}.
\end{equation}
Therefore the total bound for the second moment is
\[
qX^{1+o(1)}+X^2\log\log X\Big(\frac{1}{\log^{2024}X}+\frac{H^{2}}{(\log P_2)^{2\gamma}}
+\frac{1}{H}+\frac{\log P_1}{\log P_2}\Big).
\]
Taking
\[
P_1=\exp(\log^{\beta+\varepsilon} X), \,
P_2=\exp{\Big(\frac{\log X}{\log\log X}\Big)}, \,
H=\log^{\frac{2}{3}\gamma}X.
\]
Thus the total bound becomes
$$qX^{1+o(1)}+\frac{X^2}{\log^{\frac{2}{3}\gamma-o(1)} X}+\frac{X^2}{\log^{1-\beta-o(1)} X}.$$
\eqref{eqn: q_2m} is followed upon combining the these bounds together.

As for \eqref{eqn: Q_2m}, repeating the progress of proving \eqref{eqn: q_2m} by using \eqref{eqn: Q-ls} instead of \eqref{eqn: q-ls},
we can finish the proof.

\subsection{Proof of Proposition \ref{prop:2m}}\label{sec: 2ME}

By using Theorem \ref{thm: 2momentestimate}, to prove Proposition \ref{prop:2m}, it remains to \revisedtwo{verify} (\emph{C4}).
\revisedtwo{Let $\mathcal{E}_q$ be the set of all primitive Dirichlet characters modulo $q$. Using Lemma \ref{lemma: bound p-sum}, note that} if $\log X\asymp \log q$,
when $\exp(\log^{\beta+\varepsilon} X)\leq Y\leq \exp(\frac{\log X}{\log\log X}) $ and $\frac{Y}{\log X}\leq h\leq Y$,  we have
\begin{equation*}
	\sup_{\chi:\text{ primitive }\mod q}\Big|\sum_{Y\leq p\leq Y+h}f(p)\chi(p)\Big|\ll \frac{Y}{(\log Y)^A}.
\end{equation*}
This is (\emph{C4}) for $\gamma=A$ arbitrarily large but in certain range $\log X\asymp \log q$.

\section{Proof of Theorem \ref{thm:2ML} and Theorem \ref{thm:2ML2}}
Now we prove Theorem \ref{thm:2ML}, and Theorem \ref{thm:2ML2} can be proved via the same way.
\revisedtwo{We begin with the argument of approximate functional equation for $L(s, F\times\chi)$ at the central value} (see \cite[Theorem 5.3]{I-K04}).
For $L(\frac{1}{2}, F\times\chi)\neq 0$,
we have
\begin{multline}
	L\Big(\frac{1}{2}, F\times\chi \Big)=\frac{1}{2\pi i}\int_{c-i\infty}^{c+i\infty}L(s+\frac{1}{2},  F\times\chi)
	\frac{L_{\infty}(s+\frac{1}{2},  F\times\chi)}{L_{\infty}(\frac{1}{2}, F\times\chi)}e^{s^2}\frac{\dd s}{s}\\
	+\frac{\kappa}{2\pi i}\int_{c-i\infty}^{c+i\infty}L(s+\frac{1}{2},  \Tilde{F}\times\overline{\chi})
	\frac{L_{\infty}(s+\frac{1}{2},  \Tilde{F}\times\overline{\chi})}{L_{\infty}(\frac{1}{2}, \Tilde{F}\times\overline{\chi})}e^{s^2}\frac{\dd s}{s}
\end{multline}
with $c>\frac{1}{2}.$
By partial summation formula, we have
\begin{equation}
	L\Big(s+\frac{1}{2},  F\times\chi\Big)=(s+\frac{1}{2})\int_{1}^{\infty}\sum_{n\leq x}\lambda_F(n)\chi(n)\frac{\dd x}{x^{s+\frac{3}{2}}}.
\end{equation}
So we get
\begin{equation}
	L\Big(\frac{1}{2}, F\times\chi \Big)=\int_{1}^{\infty}\sum_{n\leq x}\lambda_F(n)\chi(n)g_\chi(x, q)\frac{\dd x}{x^{\frac{3}{2}}}+(\text{Dual integral}).
\end{equation}
\revisedtwo{where the weight function is}
\begin{equation}
	g_\chi(x, q)=\frac{1}{2\pi i}\int_{c-i\infty}^{c+i\infty}(s+\frac{1}{2})
	\frac{L_{\infty}(s+\frac{1}{2},  F\times\chi)}{L_{\infty}(\frac{1}{2}, F\times\chi)}e^{s^2}x^{-s}\frac{\dd s}{s}.
\end{equation}
Note that $L_{\infty}(s+\frac{1}{2},  F\times\chi)$ is analytic in the region $\min_{1\leq j\leq d}\Re(s+\mu_j)+\frac{1}{2}>0$, so we can shift the line of the integral to any $\Re(s)=c>\frac{1}{2}-\theta_F$ where $\theta_F$ is a positive constant satisfying \eqref{eqn: theta_f}. By Stirling's formula,  we have,
\begin{equation}\label{eqn: gxq-bound}
	\sup_{\chi\in\mathcal{F}_q}|g_\chi(x, q)|\ll_{F}\Big(\frac{q^{\frac{d}{2}}}{x}\Big)^{B},
	\quad \text{ for any } B>\frac{1}{2}-\theta_F.
\end{equation}
\revisedtwo{Note that}
\[
\sideset{}{}{\sum}_{\chi \in \mathcal{F}_q}\Big|L\big(\frac{1}{2}, F\times \chi \big)\Big|^2=\sideset{}{}{\sum}_{\chi \in \mathcal{F}_q^*}\Big|L\big(\frac{1}{2}, F\times \chi \big)\Big|^2,
\]
where $\mathcal{F}_q^*=\{\chi\in\mathcal{F}_q:L(\frac{1}{2},F\times\chi)\neq 0\}$. By using the basic inequality $|a+b|^2\leq 2(|a|^2+|b|^2)$,
it suffices to consider the first integral above and the second is estimated similarly.
So our goal is to show
\begin{equation*}
	\sideset{}{}{\sum}_{\chi \in \mathcal{F}_q^*}\Big|\int_{1}^{\infty}\sum_{n\leq x}\lambda_F(n)\chi(n)g_{\chi}(x, q)\frac{\dd x}{x^{\frac{3}{2}}} \Big|^2\ll_F \frac{q^{d/2}}{\log^{\eta} q}.
\end{equation*}
By using an hybrid form of Minkowski's inequality, we have
\begin{align*}
	\sideset{}{}{\sum}_{{\chi \in \mathcal{F}_q^*}}&\Big|\int_{1}^{\infty}\sum_{n\leq x}\lambda_F(n)\chi(n)g(x, q)\frac{\dd x}{x^{\frac{3}{2}}} \Big|^2\\
	&\leq
	\left(\int_{1}^{\infty}\Big(\sideset{}{}{\sum}_{{\chi \in \mathcal{F}_q}}\Big|\sum_{n\leq x}\lambda_F(n)\chi(n)g_\chi(x, q)\Big|^2\Big)^{\frac{1}{2}}\frac{\dd x}{x^{\frac{3}{2}}} \right)^2.
\end{align*}
We split $x$-integral into three parts $[1, \frac{q^{\frac{d}{2}}}{\log^A q}]$, $(\frac{q^{\frac{d}{2}}}{\log^A q}, q^A]$ and $(q^A,\infty)$
with some fixed large $A\geq \frac{4}{\theta_F}$. \revisedtwo{Using the basic inequality} $|a+b|^2\leq 2(|a|^2+|b|^2)$ and the large sieve inequality \eqref{eqn: q-ls}, the first part contributes at most
\begin{align*}
	&\left(\int_{1}^{\frac{q^{\frac{d}{2}}}{\log^A q}}\Big(\sideset{}{}{\sum}_{{\chi \in \mathcal{F}_q}}\Big|\sum_{n\leq x}\lambda_F(n)\chi(n)g_\chi(x, q)\Big|^2\Big)^{\frac{1}{2}}\frac{\dd x}{x^{\frac{3}{2}}} \right)^2\\
	&\ll \left(\int_{1}^{\frac{q^{\frac{d}{2}}}{\log^A q}}(qx+x^2)^{\frac{1}{2}}\sup_{\chi}|g_{\chi}(x, q)|\frac{\dd x}{x^{\frac{3}{2}}} \right)^2.
\end{align*}
Using the bound \eqref{eqn: gxq-bound} with $B=\frac{1}{2}-\theta_F+o(1)$ and $B=\frac{1}{2}-\frac{\theta_F}{2}$ respectively, it is bounded by
\begin{equation*}
	\begin{split}
		&\left(\int_{1}^{\frac{q^{\frac{d}{2}}}{\log^A q}}(qx)^{\frac{1}{2}}\Big(\frac{q^{\frac{d}{2}}}{x}\Big)^{\frac{1}{2}-\theta_F+o(1)}\frac{\dd x}{x^{\frac{3}{2}}} \right)^2
		+\left(\int_{1}^{\frac{q^{\frac{d}{2}}}{\log^A q}}\Big(\frac{q^{\frac{d}{2}}}{x}\Big)^{\frac{1}{2}-\frac{\theta_F}{2}}\frac{\dd x}{x^{\frac{1}{2}}} \right)^2\\
		&\quad\ll q^{\frac{d}{2}+\frac{1-2d\theta_F}{2}+o(1)}\left(\int_{1}^{\frac{q^{\frac{d}{2}}}{\log^A q}}\frac{\dd x}{x^{\frac{3}{2}-\theta_F-o(1)}} \right)^2+
		q^{\frac{d-d\theta_F}{2}}
		\left(\int_{1}^{\frac{q^{\frac{d}{2}}}{\log^A q}}\frac{\dd x}{x^{1-\frac{\theta_F}{2}}} \right)^2\\
		&\quad\ll  q^{\frac{d}{2}+\frac{1-2d\theta_F}{2}+o(1)}+\frac{q^{\frac{d}{2}}}{\log^2{q}}.
	\end{split}
\end{equation*}
\revisedtwo{Using Proposition \ref{prop:2m}}, the second part contributes
\begin{align*}
&	\left(\int_{\frac{q^{\frac{d}{2}}}{\log^A q}}^{q^A}\Big(\sideset{}{}{\sum}_{{\chi \in \mathcal{F}_q}}\Big|\sum_{n\leq x}\lambda_F(n)\chi(n)g_\chi(x, q)\Big|^2\Big)^{\frac{1}{2}}\frac{\dd x}{x^{\frac{3}{2}}} \right)^2\\
&	\ll \frac{1}{\log^{\eta} q}\left(\int_{\frac{q^{\frac{d}{2}}}{\log^A q}}^{\infty}\sup_{\chi}|g_{\chi}(x, q)|\frac{\dd x}{x^{\frac{1}{2}}} \right)^2.
\end{align*}
We divide the $x$-integral into $[\frac{q^{\frac{d}{2}}}{\log^A q}, q^{\frac{d}{2}}]$ and $(q^{\frac{d}{2}}, \infty)$. Using bound \eqref{eqn: gxq-bound} with $B=\frac{1}{2}-\frac{\theta_F}{2}$ and $B=2$ respectively. Hence the above is bounded by
\begin{equation*}
	\begin{split}
		&\frac{1}{\log^{\eta} q}\left(\int_{\frac{q^{\frac{d}{2}}}{\log^A q}}^{q^{\frac{d}{2}}}
		\Big(\frac{q^{\frac{d}{2}}}{x}\Big)^{\frac{1}{2}-\frac{\theta_F}{2}}\frac{\dd x}{x^{\frac{1}{2}}} \right)^2
		+\left(\int_{q^{\frac{d}{2}}}^{\infty}\Big(\frac{q^{\frac{d}{2}}}{x}\Big)^{2}\frac{\dd x}{x^{\frac{1}{2}}} \right)^2\\
		&\quad\ll \frac{q^{\frac{d}{2}-\frac{d\theta_F}{2}}}{\log^\eta q}
		\left(\int_{\frac{q^{\frac{d}{2}}}{\log^A q}}^{q^{\frac{d}{2}}}\frac{\dd x}{x^{1-\frac{\theta_F}{2}}} \right)^2+
		\frac{q^{2d}}{\log^\eta q}
		\left(\int_{q^{\frac{d}{2}}}^{\infty}\frac{\dd x}{x^{\frac{5}{2}}} \right)^2\ll \frac{q^{\frac{d}{2}}}{\log^\eta{q}}.
	\end{split}
\end{equation*}
For the third part, by taking $B=2$ and $A$ large,
\begin{align*}
	\left(\int_{q^A}^{\infty}\Big(\sideset{}{}{\sum}_{{\chi \in \mathcal{F}_q}}\Big|\sum_{n\leq x}\lambda_F(n)\chi(n)g_\chi(x, q)\Big|^2\Big)^{\frac{1}{2}}\frac{\dd x}{x^{\frac{3}{2}}} \right)^2
	&\ll \left(\int_{q^A}^{\infty} \frac{q^{\frac{d}{2}B}}{x^{\frac{1}{2}+B}}\dd x \right)^2\\
	&\ll q^{-100}.
\end{align*}
Theorem \ref{thm:2ML} is followed upon combining the above estimates.

\section{Proof of Theorem \ref{thm:2ML'} and Theorem \ref{thm:2ML2'}}\label{sec: a5a6 to a4}
\subsection{Proof of Theorem \ref{thm:2ML'}}
Let us start with the assumption (\emph{A5}). We take $\sigma=1-\frac{L}{\log^{\beta}(qT)}$ with constants  $L>0$ and $0<\beta<1$ which will be chosen later, then we get, for any $1\leq T\leq q$,
\[
\sideset{}{^*}{\sum}_{{\chi (\mod q)}}\mathcal{N}_F\Big(1-\frac{L}{\log^{\beta}qT}, T; \chi\Big)\ll_{F} \exp\Big(C_FL\frac{\log q}{\log^{\beta}(qT)}\Big)\log^A q.
\]
This means  $\prod_{\chi\mod q}^{*}L(s, F\times\chi)$ has at most $O_F\Big(\exp\Big(C_FL\frac{\log q}{\log^{\beta}qT}\Big)\log^A q\Big)$ number of zeros in the region
\[
\Re(s)\geq 1-\frac{L}{\log^{\beta}(qT)}, \quad |\Im s|\leq T.
\]
It also implies that for all modulo $q$ primitive Dirichlet character $\chi$ but except for $O_F\Big(\exp\Big(C_FL\frac{\log q}{\log^{\beta}qT}\Big)\log^A q\Big)$ number of them, $L(s, F\times\chi)\neq 0$ in the above region.
We denote the complementary set of the exceptional set of characters by $\mathcal{F}_q$, then
\begin{equation*}
\prod_{\chi \in \mathcal{F}_q}L(s, F\times\chi)\neq 0, \quad \text{ for }\Re(s)\geq 1-\frac{L}{\log^{\beta}(qT)}, \quad |\Im s|\leq T,
\end{equation*}
and
\begin{equation*}
\sideset{}{^*}
\sum_{\chi \mod q\atop \chi\notin \mathcal{F}_q}1\ll_F\exp\Big(C_FL\frac{\log q}{\log^{\beta}(qT)}\Big)\log^A q.
\end{equation*}
By using the subconvexity bound (\emph{A6}), we have
\begin{equation*}
\sideset{}{^*}
\sum_{\chi \mod q\atop \chi\notin \mathcal{F}_q}\Big|L(\frac{1}{2}, F\times\chi)\Big|^2
\ll_F q^{\frac{d}{2}}\exp\Big(C_FL\frac{\log q}{\log^{\beta}(qT)}-2\log^\delta q\Big)\log^A q.
\end{equation*}
Now we can choose $L=1/C_F$ and $\beta=1-\delta$ to obtain strong log-power savings  for the above summation over $\chi\notin \mathcal{F}_q$. For the complementary part $\chi\in \mathcal{F}_q$, Theorem \ref{thm:2ML} gives us the following estimate
\begin{equation*}
		\sideset{}{}{\sum}_{{\chi \in \mathcal{F}_q}}\Big|L\big(\frac{1}{2}, F\times \chi \big)\Big |^2\ll_{F} \frac{q^{\frac{d}{2}}}{\log^{\eta}q},
\end{equation*}
for any $\eta<1-\beta=\delta.$ Combining these two parts, we complete the proof of Theorem \ref{thm:2ML'}.

\subsection{Proof of Theorem \ref{thm:2ML2'}}
\revisedtwo{We proceed as in the proof of Theorem \ref{thm:2ML'}, but we will use assumption (\emph{A5}') instead. }

Set $\sigma=1-\frac{L'}{\log^{\beta}(QT)}$ with constants  $L'>0$ and $0<\beta<1$ which will be chosen later. Then we get, for $1\leq T\leq Q$,
\[
\sum_{Q\leq q\leq 2Q}\sideset{}{^*}{\sum}_{{\chi (\mod q)}}\mathcal{N}_F\Big(1-\frac{L'}{\log^{\beta}QT}, T; \chi\Big)\ll_{F} \exp\Big(C_F 'L'\frac{\log Q}{\log^{\beta}(QT)}\Big)\log^A Q.
\]
This means that $\prod_{Q\leq q\leq 2Q}\prod_{\chi\mod q}^{*}L(s, F\times\chi)$ has at most $$O_F\Big(\exp\Big(C_F 'L'\frac{\log Q}{\log^{\beta}(QT)}\Big)\log^A Q\Big)$$ number of zeros in  the region $\mathcal{D}_Q$:
\[
\Re(s)\geq1-\frac{L'}{\log^{\beta}(QT)}, \quad |\Im s|\leq T.
\]
It also implies that for all primitive Dirichlet character in the set $$\bigcup_{Q\leq q\leq 2Q}\{\chi(\mod q):\chi\ \text{is primitive}\}$$ but except for $O_F\Big(\exp\Big(C_F'L'\frac{\log Q}{\log^{\beta}(QT)}\Big)\log^A Q\Big)$ number of them, $L(s, F\times\chi)\neq 0$ in the above region. 
Define $\mathcal{F}_q$ by the primitive characters modulo $q$ such that \[
L(\sigma+it,F\times\chi)\neq 0, \quad \text{ for }\sigma>1-\frac{L'}{100\log^\beta(qT)},\ |t|\leq T.
\]
Clearly the above region is contained in $\mathcal{D}_Q$ for all $Q\leq q\leq 2Q$. Then we have
\begin{equation*}
	\sum_{Q\leq q\leq 2Q}	\sideset{}{^*}
	\sum_{\chi \mod q\atop \chi\notin \mathcal{F}_q}1\ll_F\exp\Big(C_F'L'\frac{\log Q}{\log^{\beta}(QT)}\Big)\log^A Q.
\end{equation*}
By using the subconvexity bound (\emph{A6}), we have
\begin{equation*}
	\sum_{Q\leq q\leq 2Q}	\sideset{}{^*}
	\sum_{\chi \mod q\atop \chi\notin \mathcal{F}_q}\Big|L(\frac{1}{2}, F\times\chi)\Big|^2\ll_F Q^{\frac{d}{2}}\exp\Big(C_F'L'\frac{\log Q}{\log^{\beta}(QT)}-2\log^\delta Q\Big)\log^A Q.
\end{equation*}
Now we can choose $L'=1/C_F'$ and $\beta=1-\delta$ to obtain arbitrarily log-power savings  for the above summation .

 For the complementary part $\chi\in \bigcup_{Q\leq q\leq 2Q}\mathcal{F}_q$, Theorem \ref{thm:2ML2} gives us
\begin{equation*}
	 \sum_{Q\leq q\leq 2Q}\sideset{}{}{\sum}_{{\chi \in \mathcal{F}_q}}\Big|L\big(\frac{1}{2}, F\times \chi \big)\Big |^2\ll_{F} \frac{Q^{\frac{d}{2}}}{\log^{\eta}Q},
\end{equation*}
for any $\eta<1-\beta=\delta.$ Therefore Theorem \ref{thm:2ML2'} follows.

\section{Proof of Theorem \ref{thm:pi_case}}\label{sec:pi_case}
Let $\pi$ be an irreducible automorphic representation over $\GL_d$ ($d\geq 3$) over $\mathbb{Q}$ with unitary central character and $L(s, \pi)$ be its associated $L$-function.
Let  $\chi$  be a primitive Dirichlet character with the conductor $q$. The twisted  $L$-function is defined by
\[
L(s, \pi \times \chi)=\sum_{n=1}^{\infty} \frac{\lambda_{\pi}(n) \chi(n)}{n^{s}}, \quad \Re(s)>1.
\]
Associated with $\pi\times \chi$, there is the Archimedean \revisedtwo{$L$-factor} defined as
\begin{equation}\label{eqn: L-factors}
    L(s, \pi_{\infty}\times \chi_{\infty})
    :={q_{\pi, \chi}}^{\frac{s}{2}}\pi^{-\frac{ds}{2}}\prod_{j=1}^{d}\Gamma\left(\frac{s+\mu_{j, \pi}+\delta_\chi}{2}\right),
\end{equation}
where $\mu_{j, \pi}$ are local parameters at $\infty$ and
$$
\delta_\chi=\left\{\begin{array}{ll}
	0, & \text { if } \chi(-1)=1, \\
	1, & \text { if } \chi(-1)=-1.
\end{array}\right.
$$
The corresponding complete  $L$ -function
$$
\Lambda(s, \pi \times \chi)=L(s, \pi \times \chi)L\left(s, \pi_{\infty} \times \chi_{\infty}\right)
$$
has an analytic continuation to the whole complex plane and also is entire of order $1$. It satisfies the following functional equation:
$$
\Lambda(s, \pi \times \chi)=\kappa_{\pi \times \chi} q_{\pi , \chi}^{\frac{1}{2}-s} \Lambda(1-s, \tilde{\pi} \times \bar{\chi}),
$$
where $\tilde{\pi}$ and $\bar{\chi}$ denotes the dual of $\pi$ and $\chi$ respectively,
$\kappa_{\pi \times \chi}$  has modulus one and  $q_{\pi , \chi}$  is the conductor of $ \pi \times \chi$.  In fact, due to the work of Bushnell and Henniart \cite{BH97}, the conductor  $q_{\pi , \chi}$  has the upper bound
\begin{equation}
	q_{\pi , \chi} \leqslant \frac{q_{\pi} q^{d}}{\left(q_{\pi}, q\right)}\ll_{\pi} q^d .
\end{equation}
\revisedtwo{This corresponds to \eqref{eqn: N-bound}.}
According to \cite{LRS95}, we have
\begin{equation}\label{eqn: theta_d}
    |\Re{\mu_{j, \pi}}|\leq \theta_d,
\end{equation}
where
$\theta_3=\frac{5}{14} $,
$\theta_4=\frac{9}{22} $ and
$\theta_d=\frac{1}{2}-\frac{1}{d^2+1}$ for $d\geq 5$.
This confirms the property \eqref{eqn: theta_f}.

These standard features of the automorphic $L$-function $L(s, \pi)$ and $L(s, \pi\times\chi)$ fit in our framework.
Now, we give a proof of Theorem \ref{thm:pi_case}. Using Theorem \ref{thm:2ML'} and Theorem \ref{thm:2ML2'}, under GRC, it suffices to verify (\emph{A2}), (\emph{A3}), (\emph{A5}) and (\emph{A6}) for $F=\pi$.
By \cite[Lemma 5.6]{Alex23}, we have
\begin{equation}
  \sum_{p\leq X}\frac{|\lambda_{\pi}(p)|^2}{p}= \log\log X+O_{\pi}(1).
\end{equation}
According to Remark \ref{rmk:p loglog bound}, we can prove (\emph{A2}) and (\emph{A3}).
By \cite[Corollary 1.4]{L-T19}, we have, for large $T\geq 2$,
\begin{equation}
  \sum_{q\leq T}\sideset{}{^*}{\sum}_{{\chi (\mod q)}}\mathcal{N}_\pi(\sigma, T; \chi)
  \ll_{\pi} T^{c_1d^2(1-\sigma)},
\end{equation}
where $c_1$ is a small absolute positive constant.
\revisedtwo{Let $Q\geq 2$ be large and $Q\leq q\leq 2Q$.
The above estimate} implies that for $1\leq T\ll Q$,
\begin{equation}
  \sideset{}{^*}{\sum}_{{\chi (\mod q)}}\mathcal{N}_\pi(\sigma, T; \chi)
  \ll_{\pi} (qT)^{c_1d^2(1-\sigma)}
\end{equation}
and
\begin{equation}
  \sum_{Q\leq q\leq 2Q}\sideset{}{^*}{\sum}_{{\chi (\mod q)}}\mathcal{N}_\pi(\sigma, T; \chi)
  \ll_{\pi} (QT)^{c_1d^2(1-\sigma)}.
\end{equation}
These verify (\emph{A5}) and (\emph{A6}). Thus Theorem \ref{thm:pi_case} is followed upon applying Theorem \ref{thm:2ML'} and Theorem \ref{thm:2ML2'}.

\section*{Acknowledgements}
\revised{The authors thank Professors Bingrong Huang and Jianya Liu for their valuable suggestions. We also want to thank the anonymous reviewers for their careful reading and constructive comments.}


\end{document}